\documentclass[11pt]{article}
\usepackage{latexsym}
\usepackage{amsfonts}
\usepackage{amsmath}
\usepackage{amssymb}
\usepackage{url}

\makeatletter
\@addtoreset{equation}{section}

\ifx\reset@font\undefined
     \let\reset@font=\relax
\fi

\def\section{\@startsection {section}{1}{\z@}{-3.5ex plus-1ex minus
       -.2ex}{2.3ex plus.2ex}{\reset@font\large\bf}}
\def\subsection{\@startsection{subsection}{2}{\z@}{-3.25ex plus-1ex
       minus-.2ex}{-.1em}{\reset@font\large\bf}}
\def\subsubsection{\@startsection{subsubsection}{3}{\z@}{-3.25ex plus
    -1ex minus-.2ex}{-.1em}{\reset@font\normalsize\bf}}
\makeatother

\date{} 
\title{ On measures of symmetry and floating bodies}
\author{Stanislaw J. Szarek (Cleveland \& Paris)
}
\def\sph{S^{n-1}}
\def\ra{\rightarrow}
\newcommand{\R}{{\mathbb{R}}}
\newcommand{\RP}{\R ^{+}}
\newcommand{\Rn}{\R ^n}
\newcommand{\N}{{\mathbb{N}}}

\newtheorem{fact}{Fact}
\newtheorem{thm}[fact]{Theorem}

\newtheorem{lemma}[fact]{Lemma}
\newtheorem{cor}[fact]{Corollary}

\newtheorem{claim}[fact]{Claim}

\begin{document}
\maketitle 
\begin{center} {\small{\sl Dedicated to Olek Pe{\l}czy{\'n}ski, a teacher and a friend}}
\end{center}
 \begin{abstract} We consider the following measure of symmetry of a convex 
$n$-dimensional body $K$: $\rho(K)$ is the smallest constant for which there is a 
point $x \in K$ such that for partitions of $K$ by an $n-1$-dimensional 
hyperplane passing through $x$ the ratio of the volumes of the two parts 
is $\le \rho(K)$. It is well known that $\rho(K)=1$ iff $K$ is symmetric. 
We establish a precise  upper bound on $\rho(K)$; this recovers 
a 1960 result of Gr\"unbaum.  We also provide 
a characterization of equality cases  (relevant to recent results of 
Nill and Paffenholz about toric varieties) 
and relate these questions  
to the concept of convex floating bodies. 
\end{abstract}

\noindent {\bf 1.  Introduction.} 
The structure of general $n$-dimensional 
(bounded) convex bodies is understood much less than that of 
the symmetric ones. 
For example, if we endow each of these classes with the appropriate 
version of the Banach-Mazur distance, then the asymptotic order 
(as $n \ra \infty$)
of the diameter of the resulting compactum has been known for the 
symmetric case since the 1981 seminal Gluskin's paper \cite{Gl}, while the  
corresponding problem in the general case is wide open, with the first 
non-trivial results having been obtained only in the last several 
years \cite{BLPS, R}.
Various invariants have been proposed to explain the difference 
between these two classes; see, e.g., \cite{Gr1} for an early survey 
of related work.  
One such measure of symmetry (or rather of asymmetry) was mentioned 
to the author by Olek Pe{\l}czynski around the turn of the millennium.  
Given convex body $K \subset \Rn$ and $x \in K$,
consider all partitions of $K$ 
into two parts by an $n-1$-dimensional hyperplane $H$ passing through 
$x$,  and let $\rho(K,x)$ be the largest ratio 
of volumes of the two parts.  Next, let  
$\rho(K) :=\min_{x\in K}{\rho(K,x)}$.  
Clearly  $\rho(K) = 1$ if $K$ is centrally symmetric. (The reverse 
implication is also true but nontrivial, even for the larger class 
of star-shaped  bodies; see \cite{S1} for details and \cite{Gro}  for 
additional references.)  
Olek's question was to establish a dimension-free upper bound on $\rho(K)$.  
I came up with an argument  (which also established a precise  
upper bound  on $\rho(K)$ for  $K \subset \Rn$ and characterized bodies,   
for which that upper bound -- call it $\rho_n$ --  is attained) and wrote it up some time 
afterwards, but then realized that the question had been considered and solved 
by Gr\"unbaum in 1960 \cite{Gr0} and so I abandoned the project.  
However, it transpired very recently  \cite{BB, NP} that Gr\"unbaum's result and  
a characterization of  $K \subset \Rn$ for which $\rho(K)=\rho_n$ 
were relevant to problems in toric geometry.  
Moreover, since  the question in 
\cite{Gr0} was stated slightly differently than above, it led to  
an apparently non-equivalent analysis of equality cases 
(see \cite{Gr0}, Remark 4(i), p. 1260), less suitable -- 
at least without additional work -- for the applications considered in \cite{NP}.  Accordingly, 
I am posting the manuscript (with added references and minor editorial changes).

\medskip \noindent  {\bf  2. More background and the results.} 
  The parameter $\rho( K )$ is related to another geometrical concept,
the convex floating bodies of $K$. [To give meaning to the formulae 
and to avoid artificial anomalies, it should be understood that $K$ -- and all bodies 
above and in what follows -- is convex, compact and has a nonempty interior.]
 Let $\delta \in (0, 1/2]$; 
slightly modifying the original definition from \cite{SW1}, let us denote by
$K^\delta$ the intersection of all half-spaces
whose complements contain at most the proportion $\delta$ of 
the volume of $K$.  As is easy to see, if we call $\phi(K)$
the largest number $\delta$ for which the convex 
floating body $K^\delta$ is nonempty, then $\phi(K)=(\rho(K)+1)^{-1}$. 
The existence of a universal (i.e., independent of $n$ and $K$) 
strictly positive lower bound for $\phi(K)$  (and, analogously, 
universal upper bound for $\rho( K )$) has been a part of the
folklore for some time \cite{Gi, SW2}.  
Here we prove the following ``isometric" result.
\begin{thm}
Let $n \in \N $ and let $K \subset \Rn$ be a convex, compact 
body with nonempty interior. Let $c$ be the centroid of $K$ and let 
$H$ be an $n-1$-dimensional hyperplane which passes through  $c$, 
and thus divides $K$ into two parts. 
Then the ratio of volumes of the two parts is $ \le (1+1/n)^n-1=: \rho_n$ 
(which is $< e -1 < 1.7183$).
Moreover, we have an equality iff $K$ is a ``pyramid," i.e., 
$K = {\rm conv}(\{v\} \cup B)$,  where $B$ is an $n-1$-dimensional 
convex body (the ``base") and $v$ the vertex, and $H$ is parallel 
to $B$.
\end{thm}
Similar statements about $\rho(K)$ and $\phi(K)$ are then simple 
consequences.

\begin{cor}
In the notation and under the hypotheses of Theorem 1 we have 

\smallskip \noindent {\rm (i) }  
$\rho(K) \le \rho_n$;  
moreover, $\rho(K,c) \le \rho_n$,  with equality 
iff $K$ is a pyramid and the maximizing 
hyperplane $H$ is parallel to a base of the pyramid.

\smallskip \noindent {\rm (ii) }  
$\phi(K) \ge  (1+1/n)^{-n} =:\delta_n$;  
moreover, $c \in K^{\delta_n}$ with $c \in \partial K^{\delta_n}$
iff $K$ is a pyramid.
\end{cor}
The estimates $\rho(K) \le \rho_n$ and $\phi(K) \ge \delta_n$ are, 
in general, best possible as seen from the example of a simplex.

\medskip \noindent {\bf 3. The proofs.} {\em Proof of Theorem 1 }  
Without loss of generality we may assume that 
the centroid $c$ of $K$ is at the origin.  
Let $\theta \in \sph$, $H = \theta^\perp$ and consider the  function
$f = f_\theta : \R \ra \RP$ defined by 
\begin{equation} \label{sections}
f(t) := {\rm vol_{n-1}} (K \cap (H+t\theta).
\end{equation} 
It then follows that $f$ is upper semi-continuous (since $K$ is closed) 
and supported on some bounded interval $[-a, b]$ with $a, b>0$.   (In fact  
it will follow from our arguments -- and is likely well known -- 
that the ratio $|a|/|b| \in [1/n,n]$.)  
Clearly,  ${\rm vol_{n}} (K) = \int_{-a}^b {f(t) \, dt }$ and  
\begin{equation} \label{rhoeq}
\rho(K,c) = \max_{\theta \in \sph} 
{\frac{\int_0^b {f(t) \, dt}}{\int_{-a}^0 {f(t) \, dt}}}.
\end{equation} 
Additionally,  the centroid being at the origin is equivalent to
\begin{equation} \label{centroid}
\int_{-a}^b {tf(t) \, dt } = 0 .
\end{equation}
The only other property of $f$ we shall need is that $h := f^{1/{(n-1)}}$ 
is concave on $[-a,b]$,  which is a consequence of the Brunn-Minkowski 
inequality. We note that the concavity implies continuity on $(-a,b)$
and lower semi-continuity, hence continuity on $[-a,b]$. 
The Theorem will follow easily from the following two claims.

\begin{claim}  Let $n \in \N$ and $a,b > 0$. 
For any $\theta \in \sph$ and for any continuous function 
$f : [-a,b] \ra \RP$  such that $h :=f^{1/{(n-1)}}$ is concave, 
there exists a closed convex body
$K \subset \Rn$  such that $f$ is obtained  from $K, \theta$ 
via {\rm (\ref{sections})}.  
If, additionally, {\rm (\ref{centroid})} holds, then $K$ may be chosen so that 
its centroid is at the origin. Finally, 
$f$ is affine with $f(-a)=0$ (or  $f(b)=0$) iff $K$ 
is a pyramid and $\theta $ is perpendicular to a base $B$; in that case, if
{\rm (\ref{centroid})} also holds,  the ratio from {\rm (\ref{rhoeq})} 
equals  $(1+1/n)^n-1$.
\end{claim}

\begin{claim}    Let $a, b >0$.  Among continuous functions on $[-a,b]$, 
strictly positive on $(-a,b)$,
verifying {\rm (\ref{centroid})}  and such that $h :=f^{1/{(n-1)}}$  is concave on
$[-a,b]$,  the largest value  of the ratio appearing in  {\rm (\ref{rhoeq})} 
is achieved iff $h$  is affine on $[-a,b]$ with $h(-a) = 0$.
\end{claim}

Claim 3 shows that investigating $\rho(K)$ and the ratio from 
(\ref{rhoeq}) for functions verifying our assumptions are fully equivalent. 
Its proof is based on elementary geometric considerations. To construct  
$K$ starting from $f$, we choose any $n-1$-dimensional convex body $B_0$ in $H =
\theta^\perp$ with ${\rm vol_{n-1}} (B_0)=1$ and $0 \in B_0$, and set 
$K:= \bigcup_{t \in [-a,b]} t\theta + h(t)B_0$. The ``only if" part in the 
last assertion follows from the analysis of equality cases in the 
Brunn-Minkowski inequality (it occurs ``essentially iff" the two sets 
are homothetic, see \cite{S1}, Theorem 6.1.1). The details are left to the reader.

\medskip 
The proof of Claim 4 is also elementary, but less obvious.  We will use  the
following lemma, variants of which exist in the literature. Similar arguments 
were employed (independently of this note) in \cite{F}, see also \cite{FG} for a more 
conceptualized application of closely related phenomena. 

\begin{lemma}  Let $M >0$, $m \in \R$ and $n \in \N$.  
We consider the  set $\mathcal{H}$ of functions $h : \R \ra \RP$ which verify
\newline {\rm (i) } the support of $h$ is the interval $[0,b]$ (for some $b>0$) 
\newline {\rm (ii) } $h$ is continuous and concave on $[0,b]$
\newline {\rm (iii) } $h(0)=1$ and the right derivative of $h$ at $0$ is $\leq m$
\newline {\rm (iv) } if $f := h^{n-1}$,  then  $\int_0^b {tf(t) \, dt } = M$.
\newline The set  $\mathcal{H}$ is nonempty iff  $m \ge -1/\sqrt{Mn(n+1)}$. 
In that case, set $\mu(h):=  \int_0^b {f(t) \, dt}$ for $h \in \mathcal{H}$.  
The minimal value  of $\mu(h)$ is attained iff $h$ is
affine on  its support $[0,b]$ with $f(b)=0$.   
The maximal value  of $\mu(h)$ is attained iff $h(t) = 1+mt $ on the support of $h$.
\end{lemma}

\noindent {\em Proof of Claim 4.}  
Since the ratio in (\ref{rhoeq})
doesn't  change if $f(\cdot)$ is replaced by $\alpha f(\beta\,\cdot)$,  it
is enough  to consider $f$'s whose support verifies  $[-1/n,1/n] \subset [-a,b]  \subset 
[-1,1]$ and such that $f(1)=1$ (for the first inclusion, see the comments
in  the paragraph following (\ref{sections}) and at the very end of this note).  
Concavity  of $f^{1/(n-1)}$ gives then a lower bound on $\int_{-a}^0 {f(t) \, dt}$ 
and an upper bound on $\|f\|_{\infty}$ (both dependent on $n$).  It follows that the
set of the functions $f$ in question is compact  (say,  in the $L_1$ metric, which is
relevant here) and hence that the supremum of the ratio in (\ref{rhoeq}) is attained.  
(This can also be proved in a variety of ways,  including from the John's
theorem.)  Let $f_0$ be such an extremal function (with support $[-a_0,b_0]$), we shall
show that it must be of the form  indicated in Claim 4.

Indeed, if $f_0$ was not affine on $[-a_0,0]$ with $f(-a_0)=0$,  we could apply 
the Lemma with $h(\cdot) = f_0(-\,\cdot)^{1/{(n-1)}}: [0,a_0] \ra \RP$,
$M=-\int_{-a_0}^0 {tf_0(t) \, dt}$  
and $m$ equal to the right derivative of $h$ at $0$ to obtain an 
extremal $h_1$ (supported on a possibly different interval $[0,a_1]$) for which 
$\mu(h_1) < \mu(h_0)$. Defining $f_1$ to coincide with $f_0$ on $[0,b_0]$ and 
with $h_1(-\, \cdot)^{n-1}$ on $[-a_1,0]$ we would get a function for which 
the ratio from (\ref{rhoeq}) was strictly larger than for $f_0$.  At the same time,
the conditions (iii) and (iv)  from the Lemma (together with our choices of $m, M$ 
would assure that,  respectively, $f_1^{1/{(n-1)}}$ was concave on $[-a_1,b_0]$ 
and that (\ref{centroid}) was satisfied.  This shows that $f_0(t) = (1+t/a_0)^{n-1}$ 
for $t \in [-a_0,0]$.  A similar argument applied to  
$h= f_0^{1/{(n-1)}}: [0,b_0] \ra \RP$, $m=1/a_0$ and the same $M$ shows that 
$f_0$ is affine on the entire interval $[-a_0,b_0]$,  which concludes the proof of 
the Claim. (In fact we showed directly that affine $h$'s give the extremal 
value of the ratio from (\ref{rhoeq}), so we didn't really need to know 
that the supremum was attained.) 
\hfill $\Box$

\bigskip \noindent {\em Sketch of the proof of Lemma 5.} First,  
let us point out that  
the condition $m \geq -1/\sqrt{Mn(n+1)}$ is equivalent to $\mathcal{H} \neq \emptyset$. 
Indeed, this is easily deduced from the observation that if $h$ is supported on $[0,b]$
and defined there by
$h(t) = 1-t/b$,  then we have the relationship  $b=\sqrt{Mn(n+1)}$. 
The assertion of the Lemma is then ``essentially obvious  from physical
considerations.''  To obtain maximal ``mass" 
$\mu(h)$ for a given  ``moment" $M$ (the constraint (iv)), we need to place the mass as
closely to the  axis $t=0$ as allowed by the concavity condition (ii) and by (iii).  To
minimize  the mass for a fixed moment, we need to place the mass as far from  $t=0$ as 
possible subject to (ii) and (iii).  This is easily formalized.  A very similar argument
allows also to determine the largest possible value of $b$ for a given $M$, and the
smallest possible value for $b$ (if at all possible) given $m$ and $M$, and to
subsequently deduce that the ratio of the two is at most $n$, thus implying the bounds on 
the ratio $|a|/|b|$ stated in the paragraph following (\ref{sections}) and 
mildly used in the proof of Claim 4.
\hfill $\Box$

\medskip \noindent{\small {\em Acknowledgement} \ Supported in part by grants from  
NSF~(U.S.A.) and by an International Cooperation Grant from Min. des  Aff. Etrang. (France) \& KBN (Poland).}


\bigskip \noindent \footnotesize{Equipe d'Analyse Fonctionnelle, 
Institut de Math\'ematiques de Jussieu-PRG,  B.C. 247, 
Universit\'e Pierre et Marie Curie,
75252 Paris, France\\
{\it and}\\
Department of Mathematics, 
Case Western Reserve University, 
Cleveland, Ohio 44106-7058, U.S.A.} 

\smallskip  \noindent  Email:  {\footnotesize\tt%
szarek@cwru.edu}\\[.5cm]


\small
\begin{thebibliography}{~~}

\bibitem{BLPS}  W. Banaszczyk, A. E. Litvak,  A. Pajor, and S. Szarek,
{\sl The Flatness Theorem for
Nonsymmetric Convex Bodies  via the Local Theory of Banach Spaces},  
Math. of Operation Research 24 (1999), no. 3, 728-750.

\bibitem{BB}  R. J. Berman and B. Berndtsson,
{\sl The volume of K\"ahler-Einstein Fano varieties and convex bodies.}
Eprint   \url{http://arxiv.org/abs/1204.1308}

\bibitem{F} M. Fradelizi, {\em Sectional bodies associated with a convex body.} 
Proc. Amer. Math. Soc. 128 (2000), 2735-2744. 

\bibitem{FG} M. Fradelizi and O. Guedon, {\sl The extreme points of subsets of 
$s$-concave probabilities and a geometric localization theorem.} Discrete
Comput. Geom. 31 (2004), 327-335. 

\bibitem {Gl} E. D. Gluskin, {\sl The diameter of Minkowski
   compactum roughly equals to $n$.}  Funct. Anal. Appl.,~15 (1981),
   57-58 (English translation).

\bibitem{Gi} A. Giannopoulos, {\sl personal communication.}

\bibitem{Gro} H. Groemer, {\sl Fourier series and spherical harmonics in convexity. }
Handbook of convex geometry, Vol. A, B, 1259-1295, North-Holland, Amsterdam, 1993. 

\bibitem{Gr0}  B. Gr\"unbaum, 
{\sl Partitions of mass-distributions and of convex bodies by hyperplanes}, Pacific
J. Math. 10 (1960), 1257-1261. 

\bibitem{Gr1}  B. Gr\"unbaum, 
{\sl Measures of symmetry for convex sets. }
In Convexity Proc. Symposia in Pure Math. vol 7 (1963), 233-270.
Amer. Math. Soc., Providence, RI.

\bibitem{NP} B. Nill and A. Paffenholz, 
{\sl On the equality case in Ehrhart's volume conjecture.} 
Eprint   \url{http://arxiv.org/abs/1205.1270}

\bibitem {R} M. Rudelson,  {\sl Distances between non--symmetric convex bodies and the
$MM^*$-estimate.}  Positivity 4 (2000), no. 2, 161-178.

\bibitem{S1} R. Schneider,  {\sl \"Uber eine Integralgleichung in der Theorie der
 konvexen K\"orper.} Math. Nachr. 44 (1970), 55-75.
 
\bibitem{S2} R. Schneider,  {\sl Convex bodies: the Brunn-Minkowski theory.} 
 Encyclopedia of Mathematics and its Applications, 44. Cambridge University Press, Cambridge 1993. 

\bibitem{SW1} C. Schuett and E. Werner, {\em The convex floating body. }
Math. Scand. 66 (1990), no. 2, 275-290.

\bibitem{SW2} C. Schuett and E. Werner, {\sl personal communication.}



\end{thebibliography}
\end{document}